\documentclass[12pt]{amsart}
\usepackage{fullpage,amsthm}

\def\CC{\mathbb{C}}
\def\FF{\mathbb{F}}
\def\QQ{\mathbb{Q}}
\def\RR{\mathbb{R}}
\def\ZZ{\mathbb{Z}}

\newtheorem{problem}{Problem}[section]
\newtheorem{prop}[problem]{Proposition}
\theoremstyle{definition}

\newtheorem{algorithm}[problem]{Algorithm}
\numberwithin{equation}{section}

\begin{document}

\title{Search techniques for root-unitary polynomials}
\author{Kiran S. Kedlaya}
\address{Department of Mathematics \\
Massachusetts Institute of Technology \\
77 Massachusetts Avenue \\
Cambridge, MA 02139}
\email{kedlaya@mit.edu}
\urladdr{http://math.mit.edu/~kedlaya}
\date{September 25, 2007}
\thanks{The author was supported by NSF grant DMS-0400747 (including the purchase of \texttt{dwork}), NSF CAREER grant DMS-0545904, and a Sloan Research Fellowship. Andre Wibisono's research project was supported during the spring 2006 semester by MIT's Undergraduate Research Opportunities Program.}

\begin{abstract}
We give an anecdotal discussion of the problem of searching
for polynomials with all roots on the unit circle, whose coefficients
are rational numbers subject to certain congruence conditions.
We illustrate with an example from a calculation in $p$-adic cohomology
made by Abbott, Kedlaya, and Roe, in which we recover the zeta function
of a surface over a finite field.
\end{abstract}

\maketitle

\section*{Introduction}

In this note, we give an anecdotal discussion of the problem of searching
for polynomials with roots on a prescribed circle whose coefficients
are rational numbers subject to certain congruence conditions.
We were led to this problem by the use of $p$-adic cohomology
to compute zeta functions of varieties over finite fields; in that
context, one is looking for certain Weil polynomials (monic integer polynomials
with complex roots all on a circle of radius $p^{i/2}$, for some
prime number $p$ and some nonnegative integer $i$), and the cohomology
calculation imposes congruence conditions on the coefficients.
In fact, the main purpose of this note is to show that in a particular example
from \cite{akr}, the conditions obtained from the cohomology calculation
indeed suffice to uniquely determine the zeta function being sought.
We also illustrate with a larger example provided by Alan Lauder.

\section{Definitions}

A polynomial $P(z) = \sum_{i=0}^n a_i z^i \in \CC[z]$ of degree $n$
is \emph{self-inversive} if there exists $u \in \CC$ with $|u| = 1$ such that
\begin{equation} \label{eq:self-inv}
a_i = u \overline{a_{n-i}} \qquad (i=0, \dots, n);
\end{equation}
it is equivalent to require the roots of $P$ to be invariant, as a multiset,
under inversion
through the unit circle.
It appears that the class of self-inversive polynomials first
occurs in a theorem of Cohn \cite{bonsall-marden, cohn}; 
it occurs naturally in the study of the
locations of roots of polynomials and their derivatives, as in the
Schur-Cohn-Marden method \cite[p.\ 150]{marden}.

We will call a polynomial
\emph{root-unitary} if its roots all lie on the unit circle;
this clearly implies self-inversivity.
This class of polynomials has been widely studied, but does not seem to have a 
standard name: the term ``unimodular polynomial'' refers to a polynomial 
whose coefficients lie on the unit circle, while
``unitary polynomial'' is often read as a synonym for ``monic
polynomial'' (particularly by speakers of French, in which a monic polynomial
is standardly a ``polyn\^ome unitaire'').

Let $P(z) \in \RR[z]$ be a real root-unitary polynomial;
then \eqref{eq:self-inv} must hold with either $u=+1$ or $u=-1$, in which case we
say $P$ is \emph{reciprocal} or \emph{antireciprocal}, respectively.
(The terms \emph{palindromic} and \emph{antipalindromic} are also sometimes used.)
If $\deg(P)$ is odd, then $P$ must be divisible by $z+1$ or $z-1$, depending on
whether $P$ is reciprocal or antireciprocal. If $\deg(P)$ is even and
$P$ is antireciprocal, then $P$ must be divisible by $(z+1)(z-1)$.
This allows reduction of many questions about real root-unitary polynomials
to the reciprocal case.

\section{The basic problem}

The basic problem is to identify rational polynomials
with roots on a prescribed circle,
given a few initial coefficients and a congruence condition on the remaining
coefficients. One can renormalize in order to talk about root-unitary polynomials;
as noted above, there is no real harm in only looking at reciprocal root-unitary
polynomials. In any case, here is the precisely formulated question we
will consider.

\begin{problem} \label{prob:basic1}
Fix positive integers $n,k,q$ with $n \geq k$.
Also fix positive integers $m_0, \dots, m_{2n}$ such that
$m_j$ divides $m_i$ for $0 \leq i \leq j \leq n$, and 
$m_i = m_{2n-i}$ for $0 \leq i \leq 2n$.
Given integers $a_0, \dots, a_{2n}$ with $a_i = q^{n-i} a_{2n-i}$ for
$0 \leq i \leq 2n$, and $a_{2n} \neq 0$,
find all polynomials $P(z)$ with all roots on the circle $|z| = \sqrt{q}$ of the form
\[
P(z) = \sum_{i=0}^{2n} (a_i + c_i m_i)z^i, 
\]
where the $c_i \in \ZZ$ must satisfy 
$c_i = q^{n-i} c_{2n-i}$ for $0 \leq i \leq 2n$, and $c_i = 0$ for $i \geq 2n-k$.
\end{problem}

The fact that Problem~\ref{prob:basic1} is a finite problem follows easily from
the estimates
\[
|a_i + c_i m_i| \leq \binom{2n}{i} q^{i/2} |a_{2n}| \qquad (i=0, \dots, 2n);
\]
when $n$ is small, these estimates carry most of the information from the
condition that $P(z \sqrt{q})$ must
be root-unitary. However, for $n$ large, this 
is quite far from true. Indeed, by \cite[Proposition~2.2.1]{dipippo-howe} (see also
\cite{sinclair} for a generalization in the context of Mahler measures),
the space of
monic root-unitary reciprocal polynomials of degree $2n$ has volume
\[
\frac{2^n}{n!} \prod_{j=1}^n \left( \frac{2j}{2j-1} \right)^{n+1-j}
\leq \frac{2^n}{n!} \prod_{j=1}^n 2^{n+1-j} = \frac{2^{(n^2+3n)/2}}{n!}
\]
whereas the space of monic reciprocal polynomials of degree $2n$
whose coefficient
of $z^i$ has norm $\leq \binom{2n}{i}$ for $i=0, \dots, 2n-1$ has volume
\[
\prod_{j=1}^n 2 \binom{2n}{j} =
\prod_{j=1}^n 2 \prod_{i=0}^{j-1} \frac{2n-i}{j-i} \geq 
\prod_{j=1}^n 2 \prod_{i=0}^{j-1} 2 =
2^{(n^2+3n)/2}.
\]
For $n$ large, these are wildly discrepant, so one expects the restriction
of root-unitarity to carry much more information than the simple bound on the
size of coefficients.

\section{Exhaustion over a tree}
\label{sec:exhaustion}

We now describe our basic approach to Problem~\ref{prob:basic1}, starting
with a change of variable also used in \cite{dipippo-howe}.
Define a polynomial $Q(z) \in \ZZ[z]$ of degree $n$
by the formula
\[
P(z) = z^n Q(z + q/z).
\]
Then for $i=0,\dots,n$,
the coefficients of $z^{n-i}, \dots, z^n$ of $Q$ are obtained from
$a_{2n-i}, \dots, a_{2n}$ by an invertible linear transformation over $\ZZ$.
We can thus reformulate Problem~\ref{prob:basic1} as follows.

\begin{problem} \label{prob:basic2}
Fix positive integers $n,k,q$ with $n \geq k$.
Also fix positive integers $m_0, \dots, m_{n}$ such that
$m_j$ divides $m_i$ for $0 \leq i \leq j \leq n$.
Given integers $b_0, \dots, b_{n}$ with $b_n \neq 0$,
find all polynomials $Q(z) \in \ZZ[z]$ with all roots real and lying in
the interval $[-2 \sqrt{q},2 \sqrt{q}]$, such that
\[
Q(z) = \sum_{i=0}^{n} (b_i + d_i m_i)z^i 
\]
for some $d_i \in \ZZ$ with $d_i = 0$ for $i \geq n-k$.
\end{problem}
Our approach to Problem~\ref{prob:basic2} is via enumeration of a certain rooted
tree. 
\begin{prop} \label{P:tree}
Fix notation as in Problem~\ref{prob:basic2}.
Then there exist sets $S_j \subseteq \ZZ^j$ for $j=0,\dots,n-k$ 
satisfying the following conditions.
\begin{enumerate}
\item[(a)]
The set $S_0$ consists of the empty $0$-tuple.
\item[(b)]
For $0 < j \leq n-k$, if $(d_{n-k-1}, \dots, d_{n-k-j}) \in S_j$,
then $(d_{n-k-1}, \dots, d_{n-k-j+1}) \in S_{j-1}$.
\item[(c)]
For $0 \leq j \leq n-k$, if $(d_{n-k-1}, \dots, d_{n-k-j}) \in S_j$, then
$Q_0(z) = \sum_{i=0}^n b_i z^i + \sum_{i=n-k-j}^{n-k-1} d_i m_i z^i$
has the property that $Q_0^{(n-k-j)}$ has all roots in $[-2\sqrt{q},
2\sqrt{q}]$.
\item[(d)]
Every tuple $(d_{n-k-1}, \dots, d_0) \in \ZZ^{n-k}$ such that
$\sum_{i=0}^n b_i z^i + \sum_{i=0}^{n-k-1} d_i m_i z^i$
has all roots in $[-2\sqrt{q},2\sqrt{q}]$ belongs to $S_{n-k}$.
\end{enumerate}
\end{prop}
\begin{proof}
Create $S_{n-k}$ by taking all solutions of Problem~\ref{prob:basic2}, then
let $S_j$ be the set of initial segments of length $j$ occurring among
elements of $S_{n-k}$. Property (c) holds by Rolle's theorem.
\end{proof}
We may identify a system of sets as in Proposition~\ref{P:tree}
with a rooted tree, where the children of a $j$-tuple in $S_j$
are its extensions to a $(j+1)$-tuple in $S_{j+1}$.
To solve Problem~\ref{prob:basic2} in practice, we perform a depth-first
enumeration of such a tree,
and read off the solutions of Problem~\ref{prob:basic2} as the elements of
$S_{n-k}$. To describe such a tree and its enumeration, it suffices
to describe how to compute the list of children of a given node.
(One could also perform a breadth-first exhaustion, but in practice this
seems to be inferior because of increased overhead.)

Note that if one wishes to decide as soon as possible whether the number of
solutions is 0, 1, or more than 1, it may be advantageous to visit the
children of a given node in ``inside-out order'' rather than in ascending
or descending order. For instance, if a given tuple can be extended by
5, 6, 7, 8, 9, we would visit these extensions in the order 7, 6, 8, 5, 9.

\section{First approach: root-finding}
\label{sec:first}

We now describe our first algorithmic approach to Problem~\ref{prob:basic2},
and its implementation \cite{kedlaya-wibisono}
in the case where $q=1$.
(See Section~\ref{sec:nonsquare} for comments on the
remaining cases.)
The implementation, 
based on an undergraduate research project
with Andre Wibisono, uses the computer algebra system \textit{SAGE}
\cite{sage}, and relies in particular on the following components.
\begin{itemize}
\item
We use the \textit{PARI} C library \cite{pari} for polynomial manipulations
over $\QQ$, including the function \texttt{polsturm} to count roots of
a polynomial in an interval using Sturm sequences. 
(This routine requires the polynomial to be
squarefree and nonvanishing at the left endpoint of the interval; one must
write a wrapper function to lift these restrictions.)
\item
We use the \textit{GNU Scientific Library (GSL)} \cite{gsl} for computing 
approximations of roots of polynomials over $\RR$.
\item
We use \textit{Cython} \cite{cython}
for migrating some critical routines into
compiled C code.
\end{itemize}
Additional components we tried out, but did not ultimately use,
include
\textit{Maple} (commercial, compiled), \textit{Numpy} (interpreted Python),
and \textit{Singular} (compiled, but its library \texttt{rootsur} for 
Sturm sequences is interpreted at runtime);
one side benefit of \textit{SAGE} is that it makes it easy to experiment
with many other systems in this manner.
Alan Lauder has done some
further experiments using \textit{Magma}.

In this approach, we take the tree in Proposition~\ref{P:tree}
to be as large as possible, by including
all tuples not forbidden by (c). This reduces to solving the following problem.
\begin{problem} \label{prob:basic3}
Let $R(z) \in \QQ[z]$ be a polynomial 
with positive leading coefficient, such that
$R'(z)$ has all roots real and in $[-2\sqrt{q},2\sqrt{q}]$. 
Find all $c \in \ZZ$ such that $R(z) + c$
has all roots real and in $[-2\sqrt{q},2\sqrt{q}]$.
\end{problem}

Let $x_1 \leq \cdots \leq x_{d-1}$ be the roots of
$R'$ counted with multiplicity, and put $x_0 = -2 \sqrt{q}$ and $x_d = 
2\sqrt{q}$.
For $i=0, \dots, d$, put $y_i = R(x_i)$.
Then the values of $c$ we want are those for which
\begin{align*}
y_{d-2i} + c \geq 0 &\quad (i=0,\dots,\lfloor d/2 \rfloor) \\
y_{d-1-2i} + c \leq 0 &\quad (i=0, \dots, \lfloor (d-1)/2 \rfloor).
\end{align*}
One can interpret this more geometrically by drawing the graph of $R$ over 
$[-2\sqrt{q},2\sqrt{q}]$. The
values of $c$ are the negatives of the integral $y$-values between the
highest local minimum and the lowest local maximum of $R$ (inclusive), provided that
we treat $2\sqrt{q}$ as a local maximum, and treat $-2\sqrt{q}$ 
as a local maximum or minimum
depending on whether $d$ is even or odd.

Our principal method for treating Problem~\ref{prob:basic3} is
to compute numerical approximations to the $x_i$ and 
$y_i$. We throw an exception if these approximations are not sufficiently
accurate, unless $x_i = x_{i+1}$ for some $i$; we can both detect and resolve
this case using exact arithmetic.

\begin{algorithm} \label{algo:seeker2}
Consider inputs as in
Problem~\ref{prob:basic3}, together with a positive integer $p$.
Using \textit{GSL},
compute numerical approximations $\tilde{x}_1 \leq \dots \leq \tilde{x}_{d-1}$ 
to the roots of $R'$, presumed (but not guaranteed)
correct to within $2^{-p}$.
For $i=1, \dots, d-1$, put
$r_i = \lfloor \tilde{x}_i 2^{p-1} - 1 \rfloor 2^{-p+1}$ and $s_i = 
r_i + 2^{-p+3}$; 
also put $r_0 = s_0 = -2\sqrt{q}$ and $r_d = s_d = 2\sqrt{q}$.
If any of the following
conditions occur for some $i \in \{1, \dots, d-1\}$:
\begin{itemize}
\item
$s_i \geq r_{i+1}$;
\item
$(-1)^{d-i} R''(r_i) > 0$;
\item
$R'(r_i)$ and $R'(s_i)$ have the same sign;
\end{itemize}
then abort or return
according as Algorithm~\ref{algo:prescreen} aborts or returns.
If none of the conditions occur,
put $l = -\infty$ and $u = +\infty$. For $i = d, d-1, \dots, 0$ in turn:
\begin{itemize}
\item if $r_i = s_i$, put $t = R(r_i)$;
\item
if $r_i < s_i$ and $d-i$ is even, let $t$ be the value computed by applying 
Algorithm~\ref{algo:subseeker} with $[r,s] = [r_i,s_i]$ and
$t_0 = -l$, then
replace $l$ by $\max\{-t,l\}$;
\item
if $r_i < s_i$ and $d-i$ is odd, let $t$ be the value computed by applying 
Algorithm~\ref{algo:subseeker} with $[r,s] = [r_i,s_i]$ and $t_0 = u$ after
replacing $R$ by $-R$, then replace $u$ by $\min\{t,u\}$;
\item
if now $l > u$, return the empty set.
\end{itemize}
Return the range $\ZZ \cap [l,u]$; this solves
Problem~\ref{prob:basic3} if not aborted.
\end{algorithm}
\begin{proof}
The only thing that needs to be noted here is that failure to 
invoke Algorithm~\ref{algo:prescreen}
ensures that the intervals $[r_1,s_1], \dots, [r_{d-1},s_{d-1}]$
are disjoint and contain one root of $R'$ apiece, so 
the input to Algorithm~\ref{algo:subseeker} is valid.
\end{proof}

In order to determine the roundings of the $y_i$, we use exact
arithmetic as follows.
\begin{algorithm} \label{algo:subseeker}
Let $R(z) \in \QQ[z]$ be a polynomial such that $R'$ has all roots
real and distinct. Let $r, s \in \QQ$ be such that 
$R''(r) \leq 0$, and the interval
$[r,s]$ contains a local maximum
of $R$ and no other roots of $R'$.
Let $t_0 \in \ZZ \cup \{+\infty\}$.
Compute
\begin{align*}
t &= \lfloor R(r) \rfloor \\
u &= \lfloor R(r) + (s-r)R'(r) \rfloor.
\end{align*}
If $t \geq t_0$, then return $t_0$ (this can be checked before
computing $u$). Otherwise,
while $t \neq u$, repeat the following:
for $v = \lceil \frac{t+u}{2} \rceil$, if $R-v$ has any roots in $[r,s]$
as determined by \texttt{polsturm}, then replace $t$ by $v$, otherwise
replace $u$ by $v-1$.
Return $t$; then for $x$ the unique root of
$R'$ contained in $[r,s]$, either $t \geq t_0$ or 
$t = \lfloor R(x) \rfloor$.
\end{algorithm}
\begin{proof}
Since $R'$ has all roots real and distinct, $x$ must be an isolated root of
$R'$. Since $x$ is a local maximum for $R$, $R'$ must undergo a sign crossing
at $x$ from positive to negative. Since $R'$ has no other roots in $[r,s]$,
$R'$ must be positive in $[r,x)$ and negative in $(x,s]$.

The roots of $R''$ interlace those of $R'$ by Rolle's theorem,
so in $(r,x]$ we  have either zero or one root of $R''$.
The root occurs if and only if there is a sign crossing; since
$R''(x) < 0$ and $R''(r) \leq 0$, we deduce that there is no root, and
$R''(z) < 0$ for all $z \in (r,x]$.

This implies that $R'(r) \geq R'(z)$ for $z \in [r,x]$; since
$R'(r) > 0$, 
\[
R(x) = R(r) + \int_r^x R'(z)\,dz 
\leq R(r) + (x-r) R'(r) \\
\leq R(r) + (s-r) R'(r).
\]
This yields the claim.
\end{proof}

Note that to a certain extent, taking $p$ small in Algorithm~\ref{algo:seeker2}
is beneficial to Algorithm~\ref{algo:subseeker}, 
because it keeps the heights of the rationals $r_i, s_i$ small. However, it
may happen that if $p$
is too small, then the gap between the initial values of $t$ and $u$ in
Algorithm~\ref{algo:subseeker} may be quite large, and a great deal of time may
be wasted narrowing the gap.

Recall that Algorithm~\ref{algo:seeker2} does not 
treat cases of Problem~\ref{prob:basic3} in 
which $R'$ has repeated roots, or $R'(-2\sqrt{q})R'(2\sqrt{q}) = 0$;
here is a simple treatment.
In practice, these cases seems to be exceedingly rare;
for instance, they do not occur at all in the 
example of Section~\ref{sec:example}. 

\begin{algorithm} \label{algo:prescreen}
Consider inputs as in
Problem~\ref{prob:basic3}.
Put $T = \gcd(R',  (z^2-4q)R'')$; 
if $T$ is constant, then abort.
Otherwise, let
$S_1, \dots, S_k$ denote the 
distinct irreducible factors of $\gcd(R',(z^2-4q)R'')$.
Determine
whether the quotients upon 
dividing $R$ by each $S_i$ are all equal to a single integer $-c$. If so, 
use \texttt{polsturm} 
to check
whether $R(z) + c$ has all roots real and in $[-2\sqrt{q},2\sqrt{q}]$; 
if so, return the singleton
set $\{c\}$. 
In all other cases, return the empty set.
This solves Problem~\ref{prob:basic3} if not aborted.
\end{algorithm}
\begin{proof}
Suppose that $T$ is nonconstant and $R(z) + c$ has all roots real and
in $[-2\sqrt{q},2\sqrt{q}]$.
Let $r$ be a root of $T$. If $r=-2\sqrt{q}$, then by Rolle's theorem,
$R(z) + c$ has a root less than or equal to $-2\sqrt{q}$, hence 
$R(-2\sqrt{q}) + c = 0$.
Similarly, if $r = 2\sqrt{q}$, then $R(2\sqrt{q}) + c = 0$. 
If $-2\sqrt{q} < r < 2\sqrt{q}$, then
$r$ is a root of $R''$ and so must be a multiple root of $R'$; by
Rolle's theorem, $r$ must be a root of $R(z) + c$. This proves the claim.
\end{proof}

\section{Second approach: power sums}
\label{sec:second}

Inspection of the enumeration of the maximal tree in some examples
suggests that it is rather bushy,
in the sense of having many vertices with many children but few deep
descendants. This in turn suggests that a more refined tree construction
might be able to achieve substantial runtime improvements.
Our second approach, implemented in \cite{kedlaya-code}
using \textit{SAGE} and components as in the previous section
(but again restricted to the case $q=1$),
does this; it is based on estimations of power sums,
as in the work of Boyd \cite{boyd} and subsequent authors
(most notably \cite{frs}) on searching for polynomials with
small Mahler measure. 

Given a polynomial $R(z) = \sum_{i=0}^n c_i z^i$ with $c_n \neq 0$,
with roots $r_1, \dots, r_n$, the \emph{power sums} of $R$ are defined as
\[
s_j = r_1^j + \cdots + r_n^j \qquad (j = 0,1,\dots).
\]
They are related to the coefficients of $R$ via the \emph{Newton identities}:
\[
j c_{n-j} + \sum_{i=0}^{j-1} c_{n-i} s_{j-i} = 0 \qquad (j=1, \dots, n).
\]
In particular, given $c_n$, one can recover $c_{n-1}, \dots, c_{n-j}$
from $s_1, \dots, s_j$ via an invertible linear transformation over $\QQ$.
Moreover, the $j$-th power sum of $R(z) + \sum_{i=0}^{n-j} c'_i z^i$ equals
$s_j - j c'_{n-j}/c_n$. Note that \textit{PARI} provides a routine
\texttt{polsym} to generate the power sums of a polynomial.

In this tree enumeration, we will generate some nodes which do not actually
belong to the tree, because they do not satisfy (c); hence our first step
when considering a proposed node will be to check (c) using 
\texttt{polsturm}. (Profiling data in some examples suggests that this
step is a bottleneck in the computation; some improvement
may be derived by instead using Sturm-Habicht sequences, as described in
\cite{bpr}, or perhaps even using real root isolation techniques.
We plan to investigate this further.)
If (c) is satisfied, and the
node is not at maximum depth, we enumerate its children by generating
and solving an instance of the following problem.

\begin{problem} \label{prob:basic4}
Given a polynomial $R(z) = \sum_{i=0}^n c_i z^i$ with $c_n \neq 0$,
and an integer $1 \leq j \leq n$, find $l,u \in \ZZ$ 
such that 
for any real numbers $d_{n-j}, \dots, d_0$ with $d_{n-j} \in \ZZ$
and $R(z) + \sum_{i=0}^{n-j} d_i z^i$ having roots in $[-2\sqrt{q},2\sqrt{q}]$, 
we have $d_{n-j} \in [l,u]$.
\end{problem}
Note that this problem is somewhat open-ended: if $l_i,u_i$ is a solution
of Problem~\ref{prob:basic4} for $i=1, \dots, k$, then so is 
$l,u$ for $l = \max_i \{l_i\}, u = \min_i \{u_i\}$. It thus suffices
to exhibit a list of inequalities satisfied by the coefficients
of a polynomial $R(z) = \sum_{i=0}^n c_i z^i$ with all roots in $[-2\sqrt{q},
2\sqrt{q}]$;
equivalently, we may exhibit inequalities satisfied 
by the power sums $s_i$ of $R$.
Here are some convenient ones; adding additional inequalities should provide
even better results, although at some point adding a new inequality
will eliminate so few cases that it will not be worth the time required
to check it.
(The linear programming approach in \cite{frs} may prove helpful
in finding good compact systems of inequalities.)

\begin{enumerate}
\item
For $i$ even,
\[
s_i - 4 q s_{i-2} \leq 0 \qquad (i \geq 2).
\]
\item
Let $T_i(z) = \sum_{k=0}^i t_{i,k} z^k$ be the polynomial
of degree $i$ for which $T_i(2\sqrt{q}\cos \theta) = 
2\sqrt{q} \cos i\theta$ (a rescaled Chebyshev polynomial of the first
kind); then
\begin{align*}
\left|\sum_{k=0}^i t_{i,k} s_k \right| &\leq 2n\sqrt{q} \qquad (i \geq 0) \\
\left|\sum_{k=0}^{i-2} t_{i-2,k} (s_{k+2} - 2 q s_k) \right| &\leq 4nq \sqrt{q}
\qquad (i \geq 2).
\end{align*}
\item
Put $s'_i = \sum_{k=0}^i \binom{i}{k} (2\sqrt{q})^{i-k} s_k$; then
\begin{align*}
s'_i s'_{i-2} - (s'_{i-1})^2 & \geq 0 \qquad (i \geq 2) \\
s'_i - 4 \sqrt{q} s'_{i-1} &\leq 0 \qquad (i \geq 1).
\end{align*}
\item
Put $s''_i = \sum_{k=0}^i \binom{i}{k} (2 \sqrt{q})^{i-k} (-1)^k s_k$; then
\begin{align*}
s''_i s''_{i-2} - (s''_{i-1})^2 &\geq 0 \qquad (i \geq 2) \\
s''_i - 4 \sqrt{q} s''_{i-1} &\leq 0 \qquad (i \geq 1).
\end{align*}
\end{enumerate}

\section{An example}
\label{sec:example}

Here is an example of the basic problem, excerpted from \cite[\S 4.2]{akr},
and some results obtained using the algorithms and implementations
 described above.

Consider the smooth quartic surface $X$ in the projective space over the finite field
$\FF_3$ defined by the homogeneous polynomial
\[
x^4 - xy^3 + xy^2w + xyzw + xyw^2 - xzw^2 + y^4 + y^3w - y^2zw + z^4 + w^4.
\]
(As described in \cite{akr}, this polynomial was chosen essentially at random
except for a skew towards sparseness.)
Since $X$ is a K3 surface, the Hodge diamond of $X$ is
\[
\begin{array}{ccccc}
& & 1 & & \\
& 0 & & 0 & \\
1 & & 20 & & 1 \\
& 0 & & 0 & \\
& & 1 & &
\end{array}
\]
and the Hodge polygon of primitive middle cohomology
has vertices $(0,0), (1,0),
(20, 19), (21, 21)$.
Consequently, the zeta function of $X$ has the form
\[
\zeta_X(T) = \exp \left( \sum_{n=1}^\infty \frac{T^n}{n} \#X(\FF_{3^n}) \right)
= \frac{1}{(1-T)(1-3T)(1-9T)R(T)},
\]
where $R(T) \in \ZZ[T]$ is a polynomial of degree 21 such that $R(0) = 1$,
the complex roots of $R$ lie
on the circle $|T| = 3^{-1}$, and (by an inequality of Mazur)
the Newton polygon of $R$ lies above the Hodge
polygon.
In particular, the polynomial $S(T) = 3 R(T/3)$ is root-unitary and
has integral coefficients.

Define 
\begin{align*}
S_0(T) &= 3T^{21} + 5T^{20} + 6T^{19} + 7 T^{18} + 5T^{17} + 4T^{16} + 2T^{15}
-T^{14} -3T^{13} - 5T^{12} \\
&\hspace{0.5cm} - 5T^{11} - 5T^{10} - 5T^9 - 3T^8 - T^7 + 2T^6
+ 4T^5 + 5T^4 + 7T^3 + 6T^2 + 5T + 3;
\end{align*}
one easily checks that $S_0$ is root-unitary.
By explicitly enumerating $X(\FF_{q^n})$ for $n \leq 5$, one finds
\[
S(T) \equiv S_0(T) \pmod{T^{6}};
\]
by a $3$-adic cohomology computation described in \cite{akr}, one finds
\[
S(T) \equiv S_0(T) \pmod{3^5}.
\]
Having performed these computations, one wants to verify whether these restrictions
suffice to ensure $S(T) = S_0(T)$. Moreover, one also wants to know 
to what extent they can be weakened
while still forcing $S(T) = S_0(T)$, as the enumeration and cohomology calculations
become significantly more cumbersome as the strength of their results is forced
to increase. 

The result obtained here is that 
already the conditions that $S(T) \in \ZZ[T]$, $S(T)$ is root-unitary,
 and
\[
S(T) \equiv S_0(T) \pmod{3^2 T^1}
\]
force $S = S_0$.
Note that already the congruence $S(T) \equiv S_0(T) \pmod{3}$ implies that
$S$ must be reciprocal rather than antireciprocal, so we may as well put
\[
P(T) = S(T)/(T+1), \qquad P_0(T) = S_0(T)/(T+1).
\]
Then the conditions we are interested in are that $P(T) \in \ZZ[T]$,
$P(T)$ is root-unitary and reciprocal, and
\begin{equation} \label{eq:cong}
P(T) \equiv P_0(T) \pmod{3^i T^j}
\end{equation}
for various $i,j$. The asserted result is that these conditions for
$i=2,j=1$ force $P = P_0$. (It turns out that $i=1$ does not suffice even
with $j=10$.)

We checked the sufficiency of the conditions for $i=2,3,4,5$ and $j=1,2,3,4,5$
by running the implementation \cite{kedlaya-wibisono}
on one Opteron 246 CPU (64-bit, 2 GHz) of the computer
\texttt{dwork.mit.edu}. The machine has 2GB of RAM.

The timings and sizes
of the computations for various initial constraints
are summarized in Table~\ref{table:timings} (using root-finding,
as in Section~\ref{sec:first}, with rounding precision $p=32$)
and Table~\ref{table:timings2} (using power sums, as in
Section~\ref{sec:second}).
Each entry consists of
the number of CPU seconds for the calculation, rounded up to the nearest
tenth of a second,
followed by
the number of leaves (terminal nodes) in the tree over which we exhausted.
Note the significant savings achieved by the second approach.
(We did some additional experiments combining the two approaches, but 
the power sum method by itself seemed to outperform hybrid methods.)

\begin{table}[ht]
\caption{Timings for recovery of $P$ given \eqref{eq:cong},
using root-finding.}
\label{table:timings}
\begin{center}
\begin{tabular}{|c|c|c|c|c|}
\hline
& $3^2$ & $3^3$ & $3^4$ & $3^5$ \\
\hline
$T^1$ & 564.2/1011788 & 2.2/3858 & 0.1/38 & 0.1/2 \\
$T^2$ & 267.9/501620 & 2.2/3784 & 0.1/38 & 0.1/2 \\
$T^3$ & 4.3/4714 & 0.1/63 & 0.1/6 & 0.1/1 \\
$T^4$ & 1.8/1838 & 0.1/51 & 0.1/6 & 0.1/1 \\
$T^5$ & 0.7/612 & 0.1/32 & 0.1/5 & 0.1/1 \\
\hline
\end{tabular}
\end{center}
\end{table}

\begin{table}[ht]
\caption{Timings for recovery of $P$ given \eqref{eq:cong},
using power sums.} 
\label{table:timings2}
\begin{center}
\begin{tabular}{|c|c|c|c|c|}
\hline
& $3^2$ & $3^3$ & $3^4$ & $3^5$ \\
\hline
$T^1$ & 1.9/1157 & 0.1/6 & 0.1/1 & 0.1/1 \\
$T^2$ & 0.6/347 & 0.1/2 & 0.1/1 & 0.1/1 \\
$T^3$ & 0.3/117 & 0.1/1 & 0.1/1 & 0.1/1 \\
$T^4$ & 0.2/53 & 0.1/1 & 0.1/1 & 0.1/1 \\
$T^5$ & 0.1/23 & 0.1/1 & 0.1/1 & 0.1/1 \\
\hline
\end{tabular}
\end{center}
\end{table}

\section{Another example}

The previous example shows the superiority of the power sum method over the
root-finding method. This suggests trying a larger example to test the
limits of the power sum method; here is an example provided by Alan Lauder.

The reciprocal polynomial
\begin{gather*}
P(T) = 2401T^{56} - 343T^{55} - 5439T^{54} - 1050T^{53} + 7156T^{52} + 5043T^{51} - 5829T^{50} - 7990T^{49} \\
+ 1437T^{48} + 6348T^{47} + 2115T^{46} - 332T^{45} - 1756T^{44} - 4639T^{43} - 1802T^{42} + 3938T^{41} \\
+ 4762T^{40} + 16T^{39} - 3366T^{38} - 2658T^{37} - 2051T^{36} + 1572T^{35} + 5810T^{34} + 2097T^{33} \\
- 5558T^{32} - 3955T^{31} + 2598T^{30} + 1931T^{29} - 831T^{28} + 1931T^{27} +
\cdots 
\end{gather*}
is root-unitary; it arises from a $7$-adic cohomology calculation of the
primitive
middle cohomology of an elliptic surface over $\FF_7$ with Hodge diamond
\[
\begin{array}{ccccc}
& & 1 & & \\
& 0 & & 0 & \\
4 & & 49 & & 4. \\
& 0 & & 0 & \\
& & 1 & &
\end{array}
\]

As in the previous example, we ask whether a reciprocal root-unitary
polynomial $P_0(T)$ satisfying
\begin{equation} \label{eq:cong2}
P(T) \equiv P_0(T) \pmod{7^i T^j}
\end{equation}
necessarily equals $P(T)$. In the following, each expression $(A/B)$
indicates that the indicated computation required $A$ CPU seconds
and encountered $B$ terminal nodes.

\begin{itemize}
\item
For $i = 2$, $P_0 = P$ is not forced for $j=28$ (0.7/15). 
\item
For $i=3$, $P_0 = P$ is forced for $j = 25$ (336.2/355435)
but not for $j=24$ (711.7/755544).
\item
For $i=4$, $P_0 = P$ is forced for 
$j=16$ (331304.0/196405710).
We were unable to find any value of $j$ for which $P_0 = P$ is not forced.
(For comparison, the complexities for $j=17,18,19,20$ were
61787.7/36665858, 12464.5/7334642, 2275.5/1349860, 392.6/232783.)
\item
For $i=5$, $P_0 = P$ is forced for $j=1$ (93.7/13513).
\end{itemize}

One can explain this behavior heuristically by imposing only the
condition that 
(in the notation of Section~\ref{sec:second})
$|\sum_{k=0}^i t_{i,k} s_k| \leq 2n$ for $i \geq 0$. 
This restriction constrains the coefficient of $T^j$ to a range of
size 
$(4 \cdot 28)/(7^4 \cdot j)$. Once $i$ is big enough that this range typically
takes only a few elements, we can expect to be able to force $P_0 = P$.

\section{The case of nonsquare $q$}
\label{sec:nonsquare}

As noted earlier, our implementations so far have only covered the
case $q=1$. It is easy to reduce to this case from any case in which
$q$ is a square. For $q$ not a square, there are several ways to proceed;
we do not know which of these is best.
\begin{itemize}
\item
One may repeat the methods as written above, but using exact arithmetic
in the quadratic field $\QQ(\sqrt{q})$.
\item
One may replace $\sqrt{q}$ by an upper approximation by a rational number
$s$ and look for polynomials with roots in $[-2s, 2s]$, then
screen out those which do not have roots in $[-2\sqrt{q},2\sqrt{q}]$.
\item
One may consider the polynomial $S$ defined by $S(z^2) = R(z) R(-z)$.
\end{itemize}

\section{Further comments}

Note that \textit{SAGE} runs primarily in the interpreted language Python,
although many of its components either are compiled libraries,
or have been migrated to C using \textit{Cython}.
(Indeed, the latter progress has been ongoing, 
and this can be detected in the runtimes
of our algorithms under different versions of \textit{SAGE}.)
It is thus fair to ask whether some additional optimization could be achieved
by porting everything to a compiled environment.
We have already built in some savings by performing most polynomial
manipulation in \textit{PARI} with limited conversions to/from \textit{SAGE},
and by porting some key subroutines into \textit{Cython}; it is not clear
how much more room there is for improvement on this front.

Our depth-first search is implemented using a queue rather than recursion;
this has the advantage of making it easily amenable to parallelization.
Although quite sophisticated strategies have been devised for 
scheduling in the context of tree traversal (e.g., \cite{zhou}),
even implementing some simple scheduling mechanisms, such as work-stealing,
would be helpful in a multiprocessor environment.
Starting with version 2.0, \textit{SAGE} includes a subsystem called
\textit{DSAGE} (Distributed \textit{SAGE}), which may facilitate this
sort of simple parallelization.

One can use similar search techniques for polynomials with roots bounded
in a convex subset of the complex plane, since the Gauss-Lucas theorem asserts
that this property is also preserved by taking derivatives. 
We have not experimented with this in any detail.

\section*{Acknowledgments}
Thanks to Alan Lauder and Chris Davis
for feedback on early versions of this paper,
and to Josh Kantor, William Stein, and Carl Witty for implementation advice.

\end{document}